\documentclass[12pt]{amsart}
\usepackage{amsmath}
\usepackage{amsthm}
\usepackage{amsfonts}
\usepackage{amssymb}
\usepackage{psfig}
\input{xy}
\xyoption{all}
\usepackage{amscd}

\theoremstyle{plain}
\newtheorem{lemma}{Lemma}
\newtheorem{cor}{Corollary}
\newtheorem{thm}{Theorem}
\newtheorem{prop}{Proposition}

\theoremstyle{definition}
\newtheorem{definition}{Definition}

\theoremstyle{remark}

\newcommand{\xym}{\ensuremath \xymatrix@1}

\newcommand{\Q}{\ensuremath \mathbb{Q}}

\newcommand{\Z}{\ensuremath \mathbb{Z}}

\newcommand{\N}{\ensuremath \mathbb{N}}

\newcommand{\isom}{\ensuremath \cong}

\begin{document}

\title[Countable Exchange and Full Exchange Rings]{Countable Exchange and Full Exchange Rings}
\author{Pace P. Nielsen}
\address{Department of Mathematics, University of California, Berkeley, CA 94720}
\email{pace@math.berkeley.edu}

\begin{abstract}
We show that a suitable ring with a ``nice'' topology, in which
convergent limits of units are units, is an $\aleph_{0}$-exchange
ring.  We generalize the argument to show that a semi-regular
ring, $R$, with a ``nice'' topology, is a full exchange ring.
Putting these results in the language of modules, we show that a
cohopfian module with finite exchange has countable exchange.
Also, all modules with Dedekind-finite, semi-regular endomorphism
rings are full exchange modules.
\end{abstract}

\maketitle

\section*{\S 1. Introduction}\label{Section:Introduction}

The exchange property for modules was first studied in 1964 by
Crawley and J\'{o}nsson \cite{CJ}, and is defined as follows.  A
right $k$-module $M_{k}$ has the {\it $\aleph$-exchange property}
if, whenever $A=M\oplus N=\bigoplus_{i\in I}A_{i}$, with
$|I|\leqslant \aleph$, then there are submodules $A_{i}'\subseteq
A_{i}$, with $A=M\oplus\left(\bigoplus_{i\in I}A_{i}' \right)$.
If $M$ has $\aleph$-exchange for all cardinals $\aleph$ then we
say $M$ has {\it full exchange}.  If the same holds just for the
finite cardinals, we say $M$ has {\it finite exchange}.  It is
easy to show that 2-exchange is equivalent to finite exchange. An
outstanding question in module theory is whether or not finite
exchange further implies full exchange.

It turns out that the finite exchange property is an endomorphism
ring invariant; putting $E=End(M_{k})$, then $M_{k}$ has finite
exchange if and only if $E_{E}$ has finite exchange.  A ring,
$R$, such that $R_{R}$ has finite exchange is called an {\it
exchange ring}, following \cite{Wa}, and this turns out to be a
left-right symmetric condition.  Nicholson \cite{N} calls a ring
{\it suitable} if, given an equation $x+y=1$, there are orthogonal
idempotents $e\in Rx$ and $f\in Ry$ with $e+f=1$.  This turns out
to be equivalent to $R$ being an exchange ring.  It is easy to
show that semi-$\pi$-regular rings\footnote{$R/J(R)$ is
$\pi$-regular, and idempotents lift modulo $J(R)$.} are suitable,
and while this is a large class it does not exhaust all exchange
rings. Any corner ring in a suitable ring is suitable, and any
direct product of suitable rings is suitable.

Continuous modules, and hence (quasi-)injective modules, always
claim the exchange property \cite{MM2}.  Further, quasi-continuous
modules with finite exchange have full exchange \cite{OR},
\cite{MM3}. There are many other classes of modules for which
finite exchange implies full exchange, including modules which are
direct sums of indecomposables \cite{ZZ}, and modules with abelian
endomorphism rings \cite{Ni}.  It also turns out that square-free
modules\footnote{No submodule is isomorphic to a square $X\oplus
X$.} with finite exchange have countable exchange \cite{MM1}.

Every endomorphism ring, $E$, is endowed with a topology, called
the {\it finite topology}, in which a basis of neighborhoods of
zero is given by annihilators of finite subsets of $M$.  One says
that a collection $\{x_{i}\}_{i\in I}\subseteq E$ of
endomorphisms is {\it summable}, if for each $m\in M$ the set
$\{i|x_{i}(m)\neq 0\}$ is finite. One may then easily define
$\sum_{i\in I}x_{i}$ as the map $m\mapsto \sum_{i\in I}x_{i}(m)$.
Central to the study of exchange modules is the following
proposition:

\begin{prop} \label{Proposition:1}
The following are equivalent:

\noindent (1) $M$ has the $\aleph$-exchange property.
\newline (2) If we have
\[
A=M\oplus N=\bigoplus_{i\in I}A_{i}
\]
with $A_{i}\cong M$ for all $i\in I$, and $|I|\leqslant \aleph$,
then there are submodules $A_{i}'\subseteq A_{i}$ such that
\[
A=M\oplus \bigoplus_{i\in I}A_{i}'.
\]
\newline (3) Given a summable family $\{x_{i}\}_{i\in I}$ of elements of
$E$, with $\sum_{i\in I}x_{i}=1$, and with $|I|\leqslant \aleph$,
then there are orthogonal idempotents $e_{i}\in Ex_{i}$ with
$\sum_{i\in I}e_{i}=1$.
\end{prop}
\begin{proof}
This is \cite[Proposition 3]{ZZ}.
\end{proof}

Now, let $R$ be a topological ring with a linear, Hausdorff
topology. This means that there is a ring topology with a basis
of zero, say $\mathfrak{U}$, consisting of left ideals, with
$\bigcap_{U\in\mathfrak{U}}U= (0)$.  We say that a collection
$\{x_{i}\}_{i\in I}\subseteq R$ is {\it summable to $r\in R$} if
there is a finite set $F'\subseteq I$ such that $\sum_{i\in
F}a_{i}-r\in U$ for all finite sets $F\supseteq F'$.  The finite
topology on $E$ is linear and Hausdorff, and this new notion of
summability agrees with the one defined above.  Following
\cite{MM1}, we can now extract from
Proposition~\ref{Proposition:1} property (3) a ring theoretic
version of $\aleph$-exchange.

\begin{definition}
Let $R$ be a ring with a linear, Hausdorff topology.  We say that
$R$ is an $\aleph${\it -exchange ring} if, given a summable
family $\{x_{i}\}_{i\in I}\subseteq R$ with $\sum_{i\in
I}x_{i}=1$, then there are summable, orthogonal idempotents
$\{e_{i}\}_{i\in I}$ with $e_{i}\in Rx_{i}$ and $\sum_{i\in
I}e_{i}=1$.\footnote{This definition differs from the one given
in [MM$_{1}$], where the word ``orthogonal'' is missing, and the
word ``complete'' is added.  From a personal correspondence with
the first author of that paper, it was made clear that the
definition given here is the one they intended.} If this holds
for all cardinalities $\aleph$, we say the ring is a {\it full
exchange ring}.
\end{definition}

Notice, a module has $\aleph$-exchange if and only if $E$ (with
the finite topology) is an $\aleph$-exchange ring.  Also notice,
in the definition above we require $\{e_{i}\}_{i\in I}$ to be a
summable family.  When trying to verify that a ring is an
$\aleph$-exchange ring, we often need to assume some condition
which forces families of this sort to be summable.  The following
is such a condition: We say a summable family $\{x_{i}\}$ is {\it
left multiple summable} if, given an arbitrary family
$\{r_{i}\}_{i\in I}$, then the collection $\{r_{i}x_{i}\}_{i\in
I}$ is also summable. We say that a topology is {\it left multiple
summable} if all summable families are left multiple summable.
Finally, we say that a topological ring, $R$, has a {\it nice
topology} if the topology is linear, Hausdorff, and left multiple
summable. One can easily show that a complete, linear, Hausdorff
topology is nice.

In this paper, we show that a suitable ring with a nice topology,
in which convergent limits of units are units, is an
$\aleph_{0}$-exchange ring. Generalizing the proof, we then show
that Dedekind-finite, regular rings with nice topologies are full
exchange rings.  We generalize the proof further to show that
$\pi$-regular, nice topological rings are full exchange rings, if
the right regular module $R_{R}$ satisfies the $(C_{2})$
property. Further, we push these arguments through the radical.
We finish by reinterpreting these results in module-theoretic
language.

\section*{\S 2. Tools for Exchange Rings}\label{Section:Tools}

Throughout this paper we let $k$ be a ring, we let $M_{k}$ be a
right $k$-module, and put $E=End(M_{k})$, which acts on the left
of $M$. All other modules will also be right $k$-modules.  If we
have two modules $N$ and $N'$ we write $N\subseteq^{\oplus}N'$ to
mean that $N$ is a direct summand of $N'$. Also throughout, we let
$R$ be a ring, $U(R)$ the group of units, and $J(R)$ the Jacobson
radical. Rings are associative with 1, and modules are unital.

In our study of $\aleph$-exchange rings, we first investigate the
behavior of idempotents in suitable rings.  To begin, we define a
useful equivalence relation on idempotents.

\begin{definition}
Let $e,e'\in R$ be idempotents.  We say that $e$ and $e'$ are
{\it left strongly isomorphic} if $e'e=e'$ and $ee'=e$.  We write
this relation as $e\sim e'$, and it is easy to check that this is
an equivalence relation.  One also has the dual notion of right
strongly isomorphic idempotents, which we denote by $e\backsim
e'$.
\end{definition}

\begin{lemma}\label{Lemma:1}
Let $e$ and $e'$ be idempotents in a ring $R$.  The following are
equivalent:

\noindent(1) $e\sim e'$.
\newline (2) $Re=Re'$.
\newline (3) $e'=e+(1-e)re$ for some $r\in R$.
\newline (4) $e'=ue$ for some $u\in U(R)$.
\newline (5) $(1-e)\backsim(1-e')$.

Furthermore, if $R=End(M_{k})$ for some module $M_{k}$, then the
following properties are also equivalent to the ones above:

\noindent (6) $\ker(e)=\ker(e')$.
\newline  (7) $(1-e)M=(1-e')M$.
\end{lemma}
\begin{proof}
The equivalence of properties $(1)$ through $(5)$ is a simple
exercise [La$_{2}$, Exercise 21.4]. $(6)\Leftrightarrow(7)$ is
easy, as is $(1)\Leftrightarrow (6)$.
\end{proof}

In the literature, two idempotents $e,e'$ are said to be {\it
isomorphic} if $eR\isom e'R$ (or equivalently, $Re\isom Re'$).
Thus, we see that if two idempotents are left (or right) strongly
isomorphic then they are isomorphic.  On the other hand, two
idempotents are both left and right strongly isomorphic if and
only if they are equal.  So, the notion of left strongly
isomorphic idempotents is a nontrivial strengthening of the
notion of isomorphic idempotents.

The equivalence in Lemma~\ref{Lemma:1} that we need the most is
$(1)\Leftrightarrow (4)$. It turns out that we can say more about
the unit in property $(4)$. In fact, by property $(3)$,
$e'=e+(1-e)re$ for some $r\in R$. Putting $u=1+(1-e)re$, we see
that $e'=ue$, and $u$ is a unit with inverse $u^{-1}=1-(1-e)re$.
Also notice, $u(1-e)=(1-e)$. So, we may strengthen property $(4)$
to read:

{\it (4$'$) $e'=ue$  for some $u\in U(R)$, with $u(1-e)=(1-e)$.}

\noindent  Throughout the rest of the paper, we will assume $(4')$
is a part of Lemma~\ref{Lemma:1}.  As an aside, although we don't
need any further properties of the unit, $u$, constructed above,
it is also true that $u(1-e')=(1-e')$, $eu=e$, $e'u=e'$, and
$(1-e)u^{-1}=1-e'$.

The next two lemmas give us computational tools we will use to
work inductively with suitable rings.

\begin{lemma}\label{Lemma:2}
Let $R$ be a suitable ring, and let $x_{1}+x_{2}+x_{3}=1$ be an
equation in $R$.  Suppose that $x_{1}$ is an idempotent. Then
there are pair-wise orthogonal idempotents $e_{1}\in Rx_{1}$,
$e_{2}\in Rx_{2}$, and $e_{3}\in Rx_{3}$, such that
$e_{1}+e_{2}+e_{3}=1$ and $x_{1}\sim e_{1}$.
\end{lemma}
\begin{proof}
Let $f=1-x_{1}$, and multiply by $f$ on the left and right of
$x_{1}+x_{2}+x_{3}=1$ to obtain $fx_{2}f+fx_{3}f=f$.  Since
corner rings in suitable rings are suitable \cite[Proposition
1.10]{N}, $fRf$ is suitable.  Hence, there are orthogonal
idempotents $f_{2}\in fRf(fx_{2}f)$ and $f_{3}\in fRf(fx_{3}f)$
summing to $f$ (the identity in $fRf$). Write
$f_{2}=fr_{2}fx_{2}f$ and $f_{3}=fr_{3}fx_{3}f$ for some
$r_{2},r_{3}\in R$.

Let $e_{2}=f_{2}r_{2}fx_{2}\in Rx_{2}$ and let
$e_{3}=f_{3}r_{3}fx_{3}\in Rx_{3}$.  By an easy calculation we
see that $e_{2}$ and $e_{3}$ are orthogonal idempotents.  Let
$e_{1}=1-e_{2}-e_{3}$, so $e_{1}$ is orthogonal to $e_{2}$ and
$e_{3}$, and we also obtain $e_{1}+e_{2}+e_{3}=1$.

We calculate
\begin{align*}
e_{1}x_{1}&=(1-e_{2}-e_{3})(1-f)=1-e_{2}-e_{3}-f+e_{2}f+e_{3}f\\
&=e_{1}-f+f_{2}+f_{3} = e_{1}-f+f=e_{1}.
\end{align*}
So $e_{1}\in Rx_{1}$.  Finally, since $fe_{2}=e_{2}$ and
$fe_{3}=e_{3}$, we see $x_{1}e_{1}=x_{1}(1-e_{2}-e_{3})=x_{1}$.
\end{proof}

\begin{lemma}\label{Lemma:3}
Let $e$ and $e'$ be idempotents in a ring $R$, with $e\sim e'$.
Assume $R$ has a linear, Hausdorff topology.  Also assume that
$e=\sum_{i\in I}g_{i}$ where $\{g_{i}\}_{i\in I}$ is a summable
family of orthogonal idempotents. Then $\{e'g_{i}\}_{i\in I}$ is
a summable family of orthogonal idempotents with $g_{i}\sim
e'g_{i}$. Further, if $e'=ue$ then $e'g_{i}=ug_{i}$. Finally, if
$f$ is any idempotent orthogonal to $e$, then $f$ is orthogonal
to each $g_{i}$.
\end{lemma}
\begin{proof}
Notice that $g_{i}e=g_{i}=eg_{i}$ and $ee'=e$.  Therefore
\[
(e'g_{i})(e'g_{j})=e'(g_{i}e)e'g_{j}=e'g_{i}(e
e')g_{j}=e'g_{i}eg_{j}=e'g_{i}g_{j}=\delta_{i,j}e'g_{i}
\]
so they are orthogonal idempotents.  Also
$g_{i}(e'g_{i})=(g_{i}e)(e'g_{i})=g_{i}eg_{i}=g_{i}$ and clearly
$(e'g_{i})g_{i}=e'g_{i}$.  Thus $g_{i}\sim e'g_{i}$.  If $e'=ue$
then $e'g_{i}=ueg_{i}=ug_{i}$.  The final statement is another
easy calculation.
\end{proof}

It will turn out that we will be working with families of
idempotents that are ``almost'' orthogonal, which we want to
modify into truly orthogonal families.  The following lemmas gives
us the mathematical framework to make this happen.

\begin{lemma}\label{Lemma:4}
Let $\{e_{i}\}_{i\in I}$ be a summable family of idempotents in a
ring $R$ with a linear, Hausdorff topology, and assume $I$ is
well-ordered. Suppose that $e_{i}e_{j}\in J(R)$ whenever $i<j$,
and that $\sum_{i\in I}e_{i}=u\in U(R)$. Then
$\{u^{-1}e_{i}\}_{i\in I}$ is a family of orthogonal idempotents,
summing to 1.
\end{lemma}
\begin{proof}
Follows from \cite[Lemma 8]{MM1}.
\end{proof}

\begin{lemma}\label{Lemma:5}
Let $\{e_{i}\}_{i\in I}$ be a summable family of idempotents in a
ring $R$ with a linear, Hausdorff topology, and assume $I$ is
well-ordered. Put $e=\sum_{i\in I}e_{i}$ and suppose that
$e_{i}e_{j}=0$ whenever $i<j$.  If $e^{n}r=0$, for some $r\in R$
and some $n\in \Z_{+}$, then we have $e_{i}r=0$ for all $i\in I$.
In particular, $er=0$.
\end{lemma}
\begin{proof}
We proceed by induction.  Since $e_{i}e_{j}=0$ for $i<j$, this
implies $e_{1}e=e_{1}$ (where 1 is the first element of $I$).
Therefore $e_{1}e^{n}=e_{1}$, and so $e_{1}r=e_{1}e^{n}r=0$. This
finishes the base case.

Now, suppose that $e_{i}r=0$ for all $i<\beta$.  Then
$er=\left(\sum_{i\geqslant \beta}e_{i}\right)r$.  Again since
$e_{i}e_{j}=0$ for $i<j$, we have
\[
e^{n-1}\left(\sum_{i\geqslant
\beta}e_{i}\right)=e^{n-2}\left(\sum_{i<
\beta}e_{i}+\sum_{i\geqslant
\beta}e_{i}\right)\left(\sum_{i\geqslant \beta}e_{i}
\right)=e^{n-2}\left(\sum_{i\geqslant
\beta}e_{i}\right)^{2}=\cdots=\left(\sum_{i\geqslant
\beta}e_{i}\right)^{n}.
\]
So,
\[
0=e_{\beta}e^{n}r=e_{\beta}e^{n-1}\left(\sum_{i\geqslant
\beta}e_{i}\right)r=e_{\beta}\left(\sum_{i\geqslant
\beta}e_{i}\right)^{n}r=e_{\beta}r.
\]
This finishes the inductive step.  It is now clear that $er=0$
also.
\end{proof}

\begin{lemma}\label{Lemma:6}
Let $R$ be an exchange ring with a linear, Hausdorff topology.
Then $J(R)$ is closed.
\end{lemma}
\begin{proof}
This is \cite[Lemma 11]{MM1}.  The lemma they prove is for
endomorphism rings, but the argument already works in this more
general situation.
\end{proof}

\begin{lemma}\label{Lemma:7}
Let $R$ be a suitable ring, and put $\overline{R}=R/J(R)$.  If
$\varepsilon\in \overline{R}\overline{x}$ is an idempotent, then
there is an idempotent $e\in Rx$ with $\overline{e}=\varepsilon$.
\end{lemma}
\begin{proof}
Follows easily from \cite[Corollary 7]{MM1}.
\end{proof}

\section*{\S 3. Countable Exchange Rings}\label{Section:CountableExchange}

The motivation for our first result comes from a simple
construction showing that $2$-exchange is equivalent to finite
exchange for modules, based upon ideas in \cite{N}. Unfortunately,
the method fails when trying to pass to countable exchange.
However, if one forces convergent limits of units to be units the
proof can be made to work as follows.

\begin{thm}\label{Theorem:1}
Let $R$ be a suitable ring with a nice topology. Also suppose that
convergent limits of units are units.  Then $R$ is an
$\aleph_{0}$-exchange ring.
\end{thm}

\begin{proof}
Let $\{x_{i}\}_{i\in \Z_{+}}$ be a summable family of elements in
$R$, with $\sum_{i=1}^{\infty}x_{i}=1$.  For notational ease, set
$y_{j}=\sum_{i>j}x_{i}$.   For each $j\in \Z_{+}$ we will
construct elements $e_{i,j}\in Rx_{i}$ (for $i\leqslant j$),
$f_{j}\in Ry_{j}$, and $v_{j}\in U(R)$ such that the following
conditions hold: (1) $\{e_{1,j},e_{2,j},\ldots, e_{j,j},f_{j}\}$
is a family of orthogonal idempotents, summing to 1, and (2)
$v_{j}e_{i,i}=e_{i,j}$ (for all $i\leqslant j$) and
$v_{j}f_{j}=f_{j}$.

Set $v_{1}=1$.  Since $R$ is suitable, the equation
$x_{1}+y_{1}=1$ implies that there are orthogonal idempotents
$e_{1,1}\in Rx_{1}$ and $f_{1}\in Ry_{1}$ with
$e_{1,1}+f_{1}=1$.  It is easy to check that condition (1) holds
for $j=1$, and condition (2) holds trivially in this case. This
finishes the base case. Suppose, by induction, we have fixed
elements $e_{i,j}\in Rx_{i}$ (for all $i\leqslant j$), $f_{j}\in
Ry_{j}$, and $v_{j}\in U(R)$ satisfying the conditions above, for
each $j\leqslant n$. Writing $f_{n}=ry_{n}$ for some $r\in R$, we
have
\[
1=e_{1,n}+\cdots + e_{n,n}+f_{n}=(e_{1,n}+\cdots +
e_{n,n})+rx_{n+1} + ry_{n+1}.
\]
Lemma~\ref{Lemma:2} allows us to pick pair-wise orthogonal
idempotents
\[
h_{1}\in R(e_{1,n}+\cdots + e_{n,n}),\qquad h_{2}\in Rrx_{n+1},
\qquad h_{3}\in Rry_{n+1}
\]
with $h_{1}+h_{2}+h_{3}=1$ and $h_{1}\sim \sum_{i=1}^{n}e_{1,n}$.
By Lemma~\ref{Lemma:1}, property ($4'$), there exists $u_{n+1}\in
U(R)$ such that $u_{n+1}(e_{1,n}+\cdots +e_{n,n})=h_{1}$ and
$u_{n+1}f_{n}=f_{n}$. Putting $e_{i,n+1}=u_{n+1}e_{i,n}\in Rx_{i}$
(for $i\leqslant n$), $e_{n+1,n+1}=h_{2}\in Rx_{n+1}$, and
$f_{n+1}=h_{3}\in Ry_{n+1}$, Lemma~\ref{Lemma:3} shows that
condition (1) above holds.

By Lemma~\ref{Lemma:1}, property (5), $(e_{n+1,n+1}+f_{n+1})$ is
right strongly isomorphic to $f_{n}$, hence
$f_{n}e_{n+1,n+1}=e_{n+1,n+1}$ and $f_{n}f_{n+1}=f_{n+1}$.
Putting $v_{n+1}=u_{n+1}v_{n}$, and remembering
$u_{n+1}f_{n}=f_{n}$, we calculate
\[
v_{n+1}f_{n+1}=(u_{n+1}v_{n})(f_{n}f_{n+1})=u_{n+1}v_{n}f_{n}f_{n+1}=
u_{n+1}f_{n}f_{n+1}=f_{n}f_{n+1}=f_{n+1}
\]
and similarly $v_{n+1}e_{n+1,n+1}=e_{n+1,n+1}$.  Finally, for
$i<n+1$,
$v_{n+1}e_{i,i}=u_{n+1}v_{n}e_{i,i}=u_{n+1}e_{i,n}=e_{i,n+1}$.
Therefore, condition (2) holds.  This finishes the inductive step.

So we have constructed elements $e_{i,j}$ (for $i\leqslant j$),
$f_{j}$, and $v_{j}$ satisfying the  properties above, for all
$j\in \Z_{+}$. Since $\{x_{i}\}_{i\in \Z_{+}}$ is summable, and
the topology is left multiple summable, the family
$\{e_{i,i}\}_{i\in\Z_{+}}$ is also summable. We put
$\varphi=\sum_{i\in \Z_{+}}e_{i,i}$. We want to prove that
$\varphi$ is a unit in $R$.

Since $\lim_{n\rightarrow \infty}y_{n}=0$, and the topology is
linear, we have $\lim_{n\rightarrow \infty}f_{n}=0$.  Therefore,
\[
\varphi=\sum_{i=1}^{\infty}e_{i,i}=\lim_{n\rightarrow
\infty}\left( \sum_{i=1}^{n}e_{i,i}+f_{n}\right)
=\lim_{n\rightarrow \infty}v_{n}^{-1}\left(\sum_{i=1}^{n}e_{i,n}
+ f_{n} \right)=\lim_{n\rightarrow \infty} v_{n}^{-1}.
\]
Convergent limits of units are units, so $\varphi$ is a unit.

Now, for $i<j$, we have
$e_{i,i}e_{j,j}=v_{j}^{-1}v_{j}e_{i,i}e_{j,j}=v_{j}^{-1}e_{i,j}e_{j,j}=0\in
J(E)$.  So, by Lemma~\ref{Lemma:4}, $\{\varphi^{-1}e_{i,i}\}_{i\in
\Z_{+}}$ is a summable, orthogonal set of idempotents, summing to
1. Finally, $\varphi^{-1}e_{i,i}\in Rx_{i}$, so $R$ satisfies the
definition of an $\aleph_{0}$-exchange ring.
\end{proof}

The converse of Theorem~\ref{Theorem:1} is not true. For example,
let $k=\Q$ and let $M_{k}=\Q^{(\N)}_{\Q}$ be the countable vector
space over $\Q$. Then $E$ is isomorphic to the ring of $\N\times
\N$ column-finite matrices over $\Q$.  One can easily construct a
limit of units in $E$ which converges to a non-unit, and yet $M$
has full exchange.

A natural question to ask is what convergent limits of units look
like in general.  We claim that in any ring with a linear,
Hausdorff topology, a convergent limit of units is always a left
non-zero-divisor.  To see this, let $w=\lim_{i\in I}w_{i}$ with
each $w_{i}$ a unit, and with $I$ well-ordered. Let $U\in
\mathfrak{U}$ be an arbitrary, open (left ideal) neighborhood of
$0$. If $wr=0$ then $\lim_{i\in I}w_{i}r=0$ and so, in
particular, for a large index $N$ we have $w_{N}r\in U$. But $U$
being a left ideal means $r=w_{N}^{-1}w_{N}r\in U$. Therefore
$r\in \bigcap_{U\in\mathfrak{U}}U=(0)$.  So $r=0$.

Theorem~\ref{Theorem:1} gives us the following chain of
corollaries.

\begin{cor}\label{Corollary:1}
Let $R$ be a suitable ring with a nice topology, and set
$\overline{R}=R/J(R)$.  If $\overline{R}_{\overline{R}}$ is
cohopfian, then $R$ is an $\aleph_{0}$-exchange ring.
\end{cor}
\begin{proof}
Let $w=\lim_{i\in I}w_{i}$, where $I$ is a well-ordered set, and
$w_{i}\in U(R)$ for each $i\in I$.  By Theorem~\ref{Theorem:1},
it suffices to show that $w$ is a unit.  An element $r\in R$ is a
unit if and only if $\overline{r}\in\overline{R}$ is a unit.
Further, by Lemma~\ref{Lemma:6}, we have $\overline{w}=\lim_{i\in
I}\overline{w}_{i}$ in the quotient topology.  Therefore it
suffices to show that $\overline{w}$ is a unit.  Since
$\overline{w}$ is a limit of units it is a left
non-zero-divisor.  By \cite[Exercise 4.16]{La2},
$\overline{R}_{\overline{R}}$ is cohopfian if and only if all
left non-zero-divisors are units.  Thus $\overline{w}$ is a unit.
\end{proof}

\begin{cor}\label{Corollary:2}
Let $R$ be a ring with a nice topology.  If $R$ is a
Dedekind-finite, semi-$\pi$-regular ring then $R$ is an
$\aleph_{0}$-exchange ring.
\end{cor}
\begin{proof}
All semi-$\pi$-regular rings are suitable rings.  So, from the
previous corollary, it suffices to show that
$\overline{R}_{\overline{R}}$ is cohopfian.

Fix $x\in \overline{R}$ which is a left non-zero-divisor.  Since
$\overline{R}$ is $\pi$-regular, fix some $n\geqslant 1$ such
that $x^{n}$ is (von Neumann) regular, say $x^{n}=x^{n}yx^{n}$
for some $y\in \overline{R}$.   Then $x^{n}(1-yx^{n})=0$.  Since
$x$ is a left non-zero-divisor so is $x^{n}$. Therefore
$1=yx^{n}$, and so $x$ is left-invertible.  From the
Dedekind-finiteness, which passes to $\overline{R}$, $x$ is
invertible.
\end{proof}

\begin{cor}\label{Corollary:3}
Let $R$ be a ring with a nice topology.  If $R$ is a strongly
$\pi$-regular ring then $R$ is an $\aleph_{0}$-exchange ring.
\end{cor}
\begin{proof}
Strongly $\pi$-regular rings are always Dedekind-finite and
$\pi$-regular.
\end{proof}

\section*{\S 4. Dedekind-finite, Regular
Rings}\label{Section:RegularRings}

When trying to push the proof of Theorem ~\ref{Theorem:1} up to
full exchange one runs into problems when passing through limit
ordinals. However, with the stronger hypothesis that $R$ is a
Dedekind-finite, regular ring, the proof goes through.

\begin{thm}\label{Theorem:2}
Let $R$ be a ring with a nice topology.  If $R$ is a
Dedekind-finite, regular ring then $R$ is a full exchange ring.
\end{thm}
\begin{proof}
Let $\{x_{i}\}_{i\in I}$ be a summable collection of
endomorphisms, summing to 1, with $I$ an indexing set of
arbitrary cardinality. Without loss of generality, we may assume
that $I$ is a well-ordered set, with first element 1, and last
element $\kappa$.  Put $y_{j}=\sum_{i>j}x_{i}$ and
$y_{j}'=y_{j}+x_{j}=\sum_{i\geqslant j}x_{i}$.

For each $j\in I$ we will inductively construct elements
$e_{i,j}\in Rx_{i}$ (for $i\leqslant j$), $f_{j}\in Ry_{j}$, and
$v_{j}\in U(R)$ such that: (1) $\{e_{i,j}\,\,\,(\forall
\,\,i\leqslant j), f_{j}\}$ is a family of orthogonal idempotents
summing to 1, and (2) $v_{j}e_{i,i}=e_{i,j}$ (for each
$i\leqslant j$) and $v_{j}f_{j}=f_{j}$.

Put $v_{1}=1$.  Since $R$ is regular it is suitable, and hence
$x_{1}+y_{1}=1$ implies that there are orthogonal idempotents
$e_{1,1}\in Rx_{1}$ and $f_{1}\in Ry_{1}$, which sum to 1.  This
completes the first step of our inductive definition. Now suppose
(by trans-finite induction) that for all $j<\alpha$ we have
constructed elements $e_{i,j}$ (for all $i\leqslant j$), $f_{j}$,
and $v_{j}$ satisfying the conditions above.  We have two cases.

{\bf Case 1.} $\alpha$ is not a limit ordinal.

In this case we proceed exactly as in the proof of Theorem 1.
Writing $f_{\alpha-1}=ry_{\alpha-1}$ for some $r\in R$, we have
\[
1=\sum_{i<\alpha}e_{i,\alpha-1}
+f_{\alpha-1}=\sum_{i<\alpha}e_{i,\alpha-1}+rx_{\alpha} +
ry_{\alpha}.
\]
Lemma~\ref{Lemma:2} allows us to pick orthogonal idempotents
\[
h_{1}\in R\left(\sum_{i<\alpha}e_{i,\alpha-1} \right),\qquad
h_{2}\in Rrx_{\alpha}, \qquad h_{3}\in Rry_{\alpha}
\]
with $h_{1}+h_{2}+h_{3}=1$ and $h_{1}\sim
\sum_{i<\alpha}e_{i,\alpha-1}$.  By Lemma~\ref{Lemma:1}, property
($4'$) , there exists $u_{\alpha}\in U(R)$ such that
$h_{1}=u_{\alpha}\left(\sum_{i<\alpha}e_{i,\alpha-1} \right)$ and
$u_{\alpha}f_{\alpha-1}=f_{\alpha-1}$. Putting
$e_{i,\alpha}=u_{\alpha}e_{i,\alpha-1}\in Rx_{i}$ (for
$i<\alpha$), $e_{\alpha,\alpha}=h_{2}\in Rx_{\alpha}$, and
$f_{\alpha}=h_{3}\in Ry_{\alpha}$, then Lemma~\ref{Lemma:3}
implies that these are orthogonal idempotents.  Also clearly
\[
\sum_{i\leqslant \alpha}e_{i,\alpha}+f_{\alpha}=1.
\]
Therefore, condition (1) holds when $j=\alpha$.  Checking that
condition (2) holds for $v_{\alpha}=u_{\alpha}v_{\alpha-1}$ is
done exactly as before. This completes the inductive definition of
the elements we need, when $\alpha$ is not a limit ordinal.

{\bf Case 2.} $\alpha$ is a limit ordinal.

This case is much harder and is where we really use the hypotheses
on $R$.  Setting $\varphi=\sum_{i<\alpha}e_{i,i}$, then since $R$
is regular there is some $\psi\in E$ with
$\varphi\psi\varphi=\varphi$, and in particular $p=1-\psi\varphi$
is an idempotent.  Putting $\varphi'=\varphi+p$, we claim that
$\varphi'$ is a unit.

First, we do a few calculations.  If $i<j<\alpha$, then
$e_{i,i}e_{j,j}=v_{j}^{-1}v_{j}e_{i,i}e_{j,j}=v_{j}^{-1}e_{i,j}e_{j,j}=0$.
Also notice that $\varphi p=0$.  So, by Lemma~\ref{Lemma:5},
$e_{i,i}p=0$ for all $i<\alpha$. Now, we show that $\varphi'$ is
a left non-zero-divisor.  To see this, suppose first that
$\varphi'\tau=0$ for some $\tau\in R$. If $\varphi\tau=0$ then
$0=\varphi'\tau=\varphi\tau + (1-\psi\varphi)\tau=\tau$.  So, we
may assume $\varphi\tau\neq 0$, and in particular there is a {\it
smallest} index $\beta$ with $e_{\beta,\beta}\tau\neq 0$.  Then
\[
0 = e_{\beta,\beta}(\varphi'\tau)=e_{\beta,\beta}\left(\sum_{i\in
[\beta,\alpha)}e_{i,i}\tau + p\tau\right) =
e_{\beta,\beta}\tau\neq 0
\]
giving a contradiction. Thus, in all cases, $\varphi'$ is a left
non-zero-divisor.  From our work in Corollary~\ref{Corollary:2} we
know that in a Dedekind-finite, regular ring any left
non-zero-divisor is a unit.  Therefore $\varphi'\in U(R)$.

For notational ease, put $v_{\alpha}'=(\varphi')^{-1}$.  From our
work above, we know that the decomposition
$\varphi'=\sum_{i<\alpha}e_{i,i} + p$ satisfies the hypotheses of
Lemma~\ref{Lemma:4}.  This yields
$\sum_{i<\alpha}v_{\alpha}'e_{i,i} + v_{\alpha}'p=1$ where the
summands are orthogonal idempotents. Put
$e_{i,\alpha}'=v_{\alpha}'e_{i,i}$ and
$f_{\alpha}'=v_{\alpha}'p$.  An easy calculation shows that
$f_{\alpha}'=p$, and in particular
$v_{\alpha}'f_{\alpha}'=f_{\alpha}'$, which we will need later.
We also claim $f_{j}f_{\alpha}'=f_{\alpha}'$ for all $j<\alpha$.
To see this we compute
\[
e_{i,j}f_{\alpha}'=v_{j}e_{i,i}f_{\alpha}'=v_{j}v_{\alpha}'^{-1}v_{\alpha}'e_{i,i}f_{\alpha}'=
v_{j}v_{\alpha}'^{-1}e_{i,\alpha}'f_{\alpha}'=0
\]
and so
\begin{equation}\label{equation:1}
f_{j}f_{\alpha}'=\left(1-\sum_{i\leqslant
j}e_{i,j}\right)f_{\alpha}'=f_{\alpha}'.
\end{equation}

Notice that we put hash marks on the idempotents we constructed.
This is because they are not quite the ones we set out to
construct. We need a few more modifications.  The first problem
with the idempotents we constructed above is that $f_{\alpha}'$
is not a left multiple of $y_{\alpha}'$.  We can fix this problem
by finding a new idempotent in $Ry_{\alpha}'$, which we will
eventually call $f_{\alpha}''$, which is right strongly
isomorphic to $f_{\alpha}'$.  The construction is as follows:

Since $R$ is regular, the principal right ideal
$y_{\alpha}'f_{\alpha}'R$ is generated by an idempotent
$g_{\alpha}$, due to \cite[Theorem 4.23]{La2}.  So there is some
$z_{\alpha}\in R$ with
$g_{\alpha}=y_{\alpha}'f_{\alpha}'z_{\alpha}$, where we may
assume $z_{\alpha}g_{\alpha}=z_{\alpha}$. Also note,
\begin{equation}\label{equation:2}
g_{\alpha}y_{\alpha}'f_{\alpha}'=y_{\alpha}'f_{\alpha}'.
\end{equation}
By definition, for $i<\alpha$ we have $f_{i}\in Ry_{i}$, and so we
can fix elements $r_{i}\in R$ with $f_{i}=r_{i}y_{i}$. For use
shortly, we also note
\begin{equation}\label{equation:3}
\lim_{i\rightarrow \alpha}y_{i}=y_{\alpha}'.
\end{equation}

Set $r_{\alpha}'=f_{\alpha}'z_{\alpha}$.  Then using
equations~\ref{equation:1} and \ref{equation:3} above, along with
left linearity, we have the following alternate definition of
$r_{\alpha}'$:
\begin{equation}\label{equation:4}
r_{\alpha}' = f_{\alpha}'z_{\alpha} =\lim_{i\rightarrow \alpha}
f_{\alpha}'z_{\alpha}= \lim_{i\rightarrow \alpha}
f_{i}f_{\alpha}'z_{\alpha} =\lim_{i\rightarrow \alpha}
r_{i}y_{i}f_{\alpha}'z_{\alpha}= \lim_{i\rightarrow \alpha} r_{i}
y_{\alpha}'f_{\alpha}'z_{\alpha}=\lim_{i\rightarrow\alpha}r_{i}g_{\alpha}.
\end{equation}
We define
$f_{\alpha}''=r_{\alpha}'y_{\alpha}'=f_{\alpha}'z_{\alpha}y_{\alpha}'$.
We first do the easy computation to show that this is an
idempotent:
\[
f_{\alpha}''f_{\alpha}''=f_{\alpha}'z_{\alpha}y_{\alpha}'f_{\alpha}'z_{\alpha}y_{\alpha}'
=f_{\alpha}'(z_{\alpha}y_{\alpha}'f_{\alpha}'z_{\alpha})y_{\alpha}'
= f_{\alpha}'(z_{\alpha}g_{\alpha})y_{\alpha}'
=f_{\alpha}'z_{\alpha}y_{\alpha}' = f_{\alpha}''.
\]
Using equations~\ref{equation:1} through \ref{equation:4} above,
we compute
\begin{align*}
f_{\alpha}''f_{\alpha}'& =r_{\alpha}'y_{\alpha}'f_{\alpha}'
=\left(\lim_{i\rightarrow\alpha}r_{i}g_{\alpha}\right)y_{\alpha}'f_{\alpha}'
=\lim_{i\rightarrow \alpha}r_{i}(g_{\alpha}y_{\alpha}'f_{\alpha}')
\\ & \! =\lim_{i\rightarrow \alpha}r_{i}y_{\alpha}'f_{\alpha}'
=\lim_{i\rightarrow \alpha}r_{i}y_{i}f_{\alpha}'
=\lim_{i\rightarrow \alpha}f_{i}f_{\alpha}' =\lim_{i\rightarrow
\alpha}f_{\alpha}' = f_{\alpha}'
\end{align*}
Also, clearly, $f_{\alpha}'f_{\alpha}''=f_{\alpha}''$.

We have shown $f_{\alpha}'\backsim f_{\alpha}''$. Therefore the
equivalence of properties (1) and (5) in Lemma~\ref{Lemma:1}
implies $(1-f_{\alpha}')\sim (1-f_{\alpha}'')$.  So, again by
Lemma~\ref{Lemma:1}, property ($4'$), pick some unit
$v_{\alpha}''$ such that
$v_{\alpha}''(1-f_{\alpha}')=1-f_{\alpha}''$ and
$v_{\alpha}''f_{\alpha}'=f_{\alpha}'$.  Set
$e_{i,\alpha}''=v_{\alpha}''e_{i,\alpha}'$, for $i<\alpha$. We
have $\sum_{i<\alpha}e_{i,\alpha}'' + f_{\alpha}''=1$, and
$\{e_{i,\alpha}''\,\,\,(\forall\,\,i<\alpha),f_{\alpha}''\}$ is a
summable family of orthogonal idempotents by Lemma~\ref{Lemma:3}.

With all the machinery we have built up, it is now an easy matter
to construct $e_{i,\alpha}$ (for all $i\leqslant \alpha$),
$f_{\alpha}$, and $v_{\alpha}$.  To do so, notice we have the
equation
\[
1=\sum_{i<\alpha}e_{i,\alpha}'' +
f_{\alpha}''=\sum_{i<\alpha}e_{i,\alpha}''+r_{\alpha}'x_{\alpha}+r_{\alpha}'y_{\alpha}.
\]
Now use exactly the same ideas as in Case 1 to construct the
elements we need.  However, there is one non-trivial step.  We
cannot put $v_{\alpha}=u_{\alpha}v_{\alpha-1}$ since $\alpha$ has
no predecessor.  Instead, we must put
$v_{\alpha}=u_{\alpha}v_{\alpha}''v_{\alpha}'$. It is clear that
$v_{\alpha}e_{i,i}=e_{i,\alpha}$ for $i<\alpha$, so we just need
to see that left multiplication by $v_{\alpha}$ acts as the
identity on $e_{\alpha,\alpha}$ and $f_{\alpha}$.  First, remember
$f_{\alpha}'=v_{\alpha}'f_{\alpha}'$.  Second, we chose
$v_{\alpha}''$ so that $v_{\alpha}''f_{\alpha}'=f_{\alpha}'$.
Third, just as in Case 1 where $u_{\alpha}$ was chosen so that
$u_{\alpha}f_{\alpha-1}=f_{\alpha-1}$, here we can choose
$u_{\alpha}$ so that $u_{\alpha}f_{\alpha}''=f_{\alpha}''$.
Finally, $e_{\alpha,\alpha}$ and $f_{\alpha}$ are both fixed by
left multiplication by $f_{\alpha}''$ and $f_{\alpha}'$ since
$(e_{\alpha,\alpha}+f_{\alpha}) \backsim f_{\alpha}''\backsim
f_{\alpha}'$.  Therefore,
\[
v_{\alpha}f_{\alpha}=(u_{\alpha}v_{\alpha}''v_{\alpha}')(f_{\alpha}'f_{\alpha})
=
u_{\alpha}(v_{\alpha}''v_{\alpha}'f_{\alpha}')f_{\alpha}=u_{\alpha}f_{\alpha}'f_{\alpha}
=
u_{\alpha}f_{\alpha}=u_{\alpha}(f_{\alpha}''f_{\alpha})=f_{\alpha}''f_{\alpha}=f_{\alpha}
\]
and similarly, $v_{\alpha}e_{\alpha,\alpha}=e_{\alpha,\alpha}$.
This finishes Case 2.

By trans-finite induction, we have constructed the elements we
wanted for all $j\in I$.  To finish the theorem, let
$e_{i}=e_{i,\kappa}$ for all $i\leqslant \kappa$. Then
$\{e_{i}\}_{i\in I}$ is a summable family of orthogonal
idempotents, summing to 1 (since $f_{\kappa}\in Ry_{\kappa}=(0)$),
with $e_{i}\in Rx_{i}$ for each $i\in I$.  This completes the
proof.
\end{proof}

\begin{cor}\label{Corollary:4}
Let $R$ be a ring with a nice topology.  If $R$ is unit-regular
then $R$ is a full exchange ring.
\end{cor}
\begin{proof}
Unit-regular rings are always regular and Dedekind-finite.
\end{proof}

We did not state Theorem~\ref{Theorem:2} in full generality so as
not to become bogged down with the details, and in an effort to
make the proof feel more natural.  Now that the basic construction
is finished we can work in a more general setting.

\begin{thm}\label{Theorem:3}
Let $R$ be a ring with a nice topology.  If $R$ is $\pi$-regular
and $R_{R}$ has $(C_{2})$ then $R$ is a full exchange ring.
\end{thm}
\begin{proof}
We need only look at how the hypothesis of regularity was used in
Theorem~\ref{Theorem:2}.  First was the fact that regularity
implied suitability.  But $R$ is suitable since $R$ is
$\pi$-regular.

Second, we needed $\varphi$ to be regular.  We know it is
$\pi$-regular, and so there is some $n\geqslant 1$, and some
$\psi\in R$, with $\varphi^{n}=\varphi^{n}\psi\varphi^{n}$.  Thus
$\varphi^{n}(1-\psi\varphi^{n})=0$.  By Lemma~\ref{Lemma:5},
$\varphi(1-\psi\varphi^{n})=0$.  In other words,
$\varphi=\varphi(\psi\varphi^{n-1})\varphi$.  Therefore,
$\varphi$ is still a regular element.

Third, we needed the fact that regularity plus
Dedekind-finiteness forces left non-zero-divisors to be units,
but this also holds in the case $R$ is $\pi$-regular.

Finally, $R$ being regular told us that $y_{\alpha}'f_{\alpha}'R$
was generated by an idempotent.  We claim that
$y_{\alpha}'f_{\alpha}'R\isom f_{\alpha}'R$, and therefore
$y_{\alpha}'f_{\alpha}'R$ will be generated by an idempotent
because of the $(C_{2})$ hypothesis.  It suffices to show that if
$y_{\alpha}'f_{\alpha}'r=0$ then $f_{\alpha}'r=0$. Using
equations~\ref{equation:1} and \ref{equation:3} above, we see
$f_{\alpha}'r=\lim_{i\rightarrow\alpha}f_{i}f_{\alpha}'r=\lim_{i\rightarrow
\alpha}
r_{i}y_{i}f_{\alpha}'r=\lim_{i\rightarrow\alpha}r_{i}y_{\alpha}f_{\alpha}'r=0.$
\end{proof}

\section*{\S 5. Lifting through the Jacobson
Radical}\label{Section:Radical}

Mohamed and M\"{u}ller have shown in \cite{MM2} that if $M$ is a
module such that $E/J(E)$ is regular and abelian, with
idempotents lifting modulo $J(E)$, then $M$ has full exchange. In
particular, they use this to establish that continuous modules
have exchange. Similarly, one way of further generalizing the
results of the previous sections is to try and lift the argument
through the Jacobson radical. The argument is actually quite easy.

\begin{thm}\label{Theorem:4}
Let $R$ be a ring with a nice topology.  Assume that $R_{R}$ has
$(C_{2})$ and $R$ is a Dedekind finite, semi-$\pi$-regular ring.
Then $R$ is a full exchange ring.
\end{thm}
\begin{proof}
First, notice that $R$ is suitable.  If one works through the
proofs of Theorems~\ref{Theorem:2} and \ref{Theorem:3}, the only
other point which needs some modification is the choice of the
idempotent $p$.  There are two properties we need $p$ to satisfy.
First, we need $\varphi p=0$, so that the calculation showing
$p=f_{\alpha}'$ will work, and also so $e_{i,i}p=0$ for all
$i<\alpha$.  Second, we need $\varphi'=\varphi+p$ to be a unit.

By Lemma~\ref{Lemma:6}, we have that
$\{\overline{e}_{i,i}\}_{i<\alpha}$ is summable in the quotient
topology of $R/J(R)$, summing to
$\overline{\widetilde{\varphi}}$.  Since $R/J(R)$ is
$\pi$-regular, the argument in Theorem~\ref{Theorem:3} shows that
$\overline{\varphi}$ is regular.  Hence, there is some $\psi\in
R$ with $\varphi-\varphi\psi\varphi\in J(R)$.  Since idempotents
lift modulo $J(R)$, and since $1-\psi\varphi$ is an idempotent
modulo $J(R)$, we can pick an idempotent $\widetilde{p}\in R$
(not quite the one we want) with
$\widetilde{p}-(1-\psi\varphi)\in J(R)$. Put
$\widetilde{\varphi}=\varphi+\widetilde{p}$.

We want to show $\widetilde{\varphi}$ is a unit in $R$, and so it
suffices to show that $\overline{\widetilde{\varphi}}$ is a unit
in $R/J(R)$.  But because of how $\widetilde{p}$ was chosen, the
same argument in Theorems~\ref{Theorem:2} and \ref{Theorem:3},
which showed $\varphi'$ was a unit, will now show that
$\overline{\widetilde{\varphi}}$ is a unit.  To make things
explicit, we will repeat the argument here.

Since $R/J(R)$ is $\pi$-regular and Dedekind-finite, it suffices
to show that $\overline{\widetilde{\varphi}}$ is a left
non-zero-divisor.  Suppose $\overline{\widetilde{\varphi}\tau}=0$
for some $\tau\in R$.  If $\overline{\varphi \tau}=0$, then since
$\widetilde{p}-(1-\psi\varphi)\in J(R)$, we have
$0=\overline{\widetilde{\varphi}\tau}=\overline{\varphi \tau +
(1-\psi\varphi)\tau}=\overline{\tau}$.  Therefore, we may assume
$\overline{\varphi\tau}\neq 0$, and in particular there is a
smallest index $\beta$, with $\overline{e_{\beta,\beta}\tau}\neq
0$.  Now, $\overline{\varphi \widetilde{p}}=0$, and so
Lemma~\ref{Lemma:5} implies that
$\overline{e_{i,i}\widetilde{p}}=0$ for all $i<\alpha$.
Therefore, working modulo $J(R)$, we calculate
\[
0 \equiv e_{\beta,\beta}(\varphi'\tau)\equiv
e_{\beta,\beta}\left(\sum_{i\in [\beta,\alpha)}e_{i,i}\tau +
p\tau\right)\equiv e_{\beta,\beta}\tau \notin 0+J(R).
\]
This contradiction shows that $\overline{\widetilde{\varphi}}$ is
a left non-zero-divisor, and hence a unit.

In our work above we found that $e_{i,i}\widetilde{p}\in J(R)$
for $i<\alpha$. Then, by Lemma~\ref{Lemma:4}, the collection
$\left\{\left(\widetilde{\varphi} \right)^{-1}
e_{i,i}\,\,\,(\forall\,\, i<\alpha), \left(\widetilde{\varphi}
\right)^{-1}\widetilde{p}\right\}$ consists of orthogonal
idempotents, summing to 1.  Put $p=\left(\widetilde{\varphi}
\right)^{-1}\widetilde{p}$. Since $\left(\widetilde{\varphi}
\right)^{-1} e_{i,i}p=0$ we have $e_{i,i}p=0$, and in particular
$\varphi p=0$.

Set $\varphi'=\varphi+p$.  Suppose that $\varphi \tau=0$ for some
$\tau\in R$. Then,
\[
\tau=\left(\sum_{i<\alpha}\left(\widetilde{\varphi} \right)^{-1}
e_{i,i} + p \right)\tau =\left(\widetilde{\varphi} \right)^{-1}
\varphi \tau + p\tau = p \tau=\varphi'\tau.
\]
Notice that we can push this equation down to $R/J(R)$.  Showing
$\varphi'=\varphi+p$ is a unit is now a simple matter by copying
the ideas used in the proof that $\widetilde{\varphi}$ is a unit.
\end{proof}

One also has another way to lift the argument through the radical.

\begin{cor}\label{Corollary:5}
Let $R$ be a ring with a nice topology.  If $R$ is a
Dedekind-finite, semi-regular ring then $R$ is a full exchange
ring.
\end{cor}
\begin{proof}
Let $\{x_{i}\}_{i\in I}$ be a summable family of idempotents,
summing to 1.  Let $I$ be well-ordered as usual. Putting
$\overline{R}=R/J(R)$, then we see by Lemma~\ref{Lemma:6} that
$\overline{R}$ is a topological ring in the quotient topology
with a linear, Hausdorff topology.  Further,
$\{\overline{x}_{i}\}_{i\in I}$ is a summable family summing to
1, and is left-multiple summable, since $\{x_{i}\}_{i\in I}$ is.
Therefore, the same argument as used in the proof of
Theorem~\ref{Theorem:2} shows that we can find orthogonal
idempotents $\varepsilon_{i}\in \overline{R}\overline{x_{i}}$
summing to 1.

By Lemma~\ref{Lemma:7}, we can lift each $\varepsilon_{i}$ to an
idempotent $e_{i}\in Rx_{i}$.  These are still summable
idempotents, summing to a unit (since, modulo $J(R)$, they sum to
1).  Letting $u=\sum_{i\in I}e_{i}$, then Lemma~\ref{Lemma:4} says
that $\{u^{-1}e_{i}\}_{i\in I}$ is a summable family of orthogonal
idempotents summing to 1.  Clearly, $u^{-1}e_{i}\in Rx_{i}$, so
we are done.
\end{proof}

Using the same ideas, we also have

\begin{cor}\label{Corollary:6}
Let $R$ be a ring with a nice topology.  If $R$ is a
Dedekind-finite, semi--$\pi$-regular ring, and
$\overline{R}_{\overline{R}}$ has $(C_{2})$, then $R$ is a full
exchange ring.
\end{cor}

\section*{\S 6. Exchange Modules}\label{Section:Modules}

What do the previous theorems say concerning finite exchange
modules?  We have the following unsettling fact, motivated by
\cite[Proposition 8.11]{La1}.

\begin{lemma}\label{Lemma:8}
Let $M_{k}$ be a module, and $E=End(M_{k})$, as usual. If $E_{E}$
is cohopfian, or respectively has $(C_{2})$, then so does $M$.
The converses do not hold.
\end{lemma}
\begin{proof}
First, suppose that $E_{E}$ is cohopfian.  Let $x\in E$ be an
injective endomorphism on $M$.  If $xr=0$ for some $r\in E$, then
$xr(m)=0$ for all $m\in M$.  But, $x$ being injective implies
$r(m)=0$ for all $m\in M$.  Therefore, $r=0$.  Since $r$ was
arbitrary, $x$ is a left non-zero-divisor.  Therefore, since
$E_{E}$ is cohopfian, $x$ is a unit.  This shows that $M$ is
cohopfian.

Now instead suppose that $E_{E}$ has $(C_{2})$.  Consider the
situation where $N'\isom N\subseteq^{\oplus}M$.  Let $e\in E$ be
an idempotent with $e(M)=N$, and let $\varphi:N\rightarrow N'$ be
an isomorphism. Without loss of generality, we may assume
$\varphi\in E$ by setting $\varphi$ equal to 0 on $(1-e)(M)$.

Consider the map, $eE\rightarrow \varphi eE$, given by left
multiplication by $\varphi$.  Clearly this is surjective. To show
injectivity, suppose that $\varphi e r=0$ for some $r\in E$. Then
$\varphi er(m)=0$ for all $m\in M$.  In particular,
$\varphi(er(M))=0$.  But $er(M)\subseteq e(M)$ and $\varphi$ is
injective on $e(M)=N$, therefore $er(M)=0$.  But then $er=0$. This
shows injectivity.

Thus $\varphi eE$ is isomorphic to $eE$, a direct summand of
$E_{E}$.  Therefore $\varphi eE$ is generated by an idempotent,
say $f$.  Clearly $f\varphi e=\varphi e$, and $f=\varphi e y$ for
some $y\in E$. So $f(M)=\varphi ey(M)\subseteq \varphi e(M)=N'$,
and $f(M)\supseteq f(\varphi e(M))=\varphi e(M)=N'$.  Therefore
$N'=f(M)$ is a direct summand.

A single counter-example will show that both converses do not
hold. Let $k=\mathbb{Z}$ and let $M$ be the Pr\"{u}fer $p$-group,
for any prime $p$. Then $E$ is isomorphic to the ring of $p$-adic
integers.  $M$ is cohopfian while $E$ is not, by
\cite[Proposition 8.11]{La1}.  Notice that the only idempotents in
$E$ are 0 and 1.  Thus, the only direct summands in either
$M_{k}$ or $E_{E}$ are the trivial ones.  One easily sees that
multiplication by $p$ yields $pE\isom E_{E}$, but $pE$ is not a
summand.  Therefore $E_{E}$ does not have the $(C_{2})$ property.
On the other hand, any submodule isomorphic to $M$ must contain
elements killed by multiplication by $p$, and hence must equal
$M$. Thus, all submodules of $M$ isomorphic to $M$ are summands,
and all submodules of $M$isomorphic to $(0)$ equal $(0)$.  Hence
$M$ has the $(C_{2})$ property.
\end{proof}

Due to this lemma, it would appear that one could not work with
the weaker notion of a cohopfian module and hope to prove a
theorem analogous to Theorem~\ref{Theorem:1}.  However, in
endomorphism rings, limits of units are very special.

\begin{thm}\label{Theorem:5}
Let $M$ be a cohopfian module with finite exchange.  Then $M$ has
countable exchange.
\end{thm}
\begin{proof}
In the endomorphism ring, $E$, a limit of units must be an
injective endomorphism (since nothing in the limit process has a
kernel).  But then the cohopfian condition forces this
endomorphism to be an isomorphism, or in other words a unit in
$E$.  Thus convergent limits of units are units.  So $M$ has
countable exchange from Theorem~\ref{Theorem:1}.
\end{proof}

Can one also tweak Theorem~\ref{Theorem:4} so we are working with
the weaker hypothesis that $M$ has the $(C_{2})$ property? The
answer is yes.

\begin{thm}\label{Theorem:6}
Let $M$ be a module with the $(C_{2})$ property, and a
Dedekind-finite, semi-$\pi$-regular endomorphism ring.  Then $M$
has full exchange.
\end{thm}
\begin{proof}
Following Theorem~\ref{Theorem:4}, with $R=E$, the only thing we
need to do differently is find an idempotent $f_{\alpha}''\in
Ey_{\alpha}'$ with $f_{\alpha}'\backsim f_{\alpha}''$.

Consider the map $y_{\alpha}':f_{\alpha}'(M)\rightarrow
y_{\alpha}'f_{\alpha}'(M)$, given by left-multiplication by
$y_{\alpha}'$.  It is clearly surjective.  We have
$f_{\alpha}'=\lim_{i\rightarrow \alpha}f_{i}f_{\alpha}'
=\lim_{i\rightarrow \alpha}r_{i}y_{i}f_{\alpha}'
=\lim_{i\rightarrow \alpha}r_{i}y_{\alpha}'f_{\alpha}'$, and so
the map above must also be injective.  From the $(C_{2})$
hypothesis, we have that
$y_{\alpha}'f_{\alpha}'(M)=g_{\alpha}(M)$ for some idempotent
$g_{\alpha}$.

Define $r_{\alpha}'$ by the rule
$r_{\alpha}'|_{(1-g_{\alpha})(M)}=0$ and
$r_{\alpha}'|_{g_{\alpha}(M)=y_{\alpha}'f_{\alpha}'(M)}=\lim_{i\rightarrow
\alpha}r_{i}$.  While it is true that $\lim_{i\rightarrow
\alpha}r_{i}$ does not necessarily converge in general, it does
converge on $y_{\alpha}'f_{\alpha}'(M)$ since
$\lim_{i\rightarrow\alpha}r_{i}y_{\alpha}'f_{\alpha}'(m)
=\lim_{i\rightarrow\alpha}r_{i}y_{i}f_{\alpha}'(m)
=\lim_{i\rightarrow\alpha}f_{i}f_{\alpha}'(m)
=f_{\alpha}'(m)$.\footnote{One should now also check that
$r_{\alpha}'$ is a well-defined homomorphism, which we leave to
the reader.}

We put $f_{\alpha}''=r_{\alpha}'y_{\alpha}'$.  We first check that
it is an idempotent.  Given $m\in M$, we can write
$y_{\alpha}'(m)=g_{\alpha}y_{\alpha}'(m)+(1-g_{\alpha})y_{\alpha}'(m)
=y_{\alpha}'f_{\alpha}'(m')+ (1-g_{\alpha})y_{\alpha}'(m)$ for
some $m'\in M$. Then,
\begin{align*}
f_{\alpha}'' f_{\alpha}''(m)&
=r_{\alpha}'y_{\alpha}'r_{\alpha}'(y_{\alpha}'f_{\alpha}'(m')+
(1-g_{\alpha})y_{\alpha}'(m))
=r_{\alpha}'y_{\alpha}'(r_{\alpha}'y_{\alpha}'f_{\alpha}')(m')\\&
=
r_{\alpha}'y_{\alpha}'\left(\lim_{i\rightarrow\alpha}r_{i}y_{i}f_{\alpha}'\right)(m')
=r_{\alpha}'y_{\alpha}'\left(\lim_{i\rightarrow\alpha}f_{i}f_{\alpha}'\right)(m')=r_{\alpha}'y_{\alpha}f_{\alpha}(m')
\\& =r_{\alpha}'(y_{\alpha}'f_{\alpha}'(m')+
(1-g_{\alpha})y_{\alpha}'(m))=r_{\alpha}'y_{\alpha}'(m)=f_{\alpha}''(m).
\end{align*}
So $f_{\alpha}''f_{\alpha}''=f_{\alpha}''$. A similar computation
shows that $f_{\alpha}''$ and $f_{\alpha}'$ are right strongly
isomorphic.  The rest of the proof follows
Theorem~\ref{Theorem:4}.
\end{proof}

Theorem~\ref{Theorem:2}, and Corollaries~\ref{Corollary:1}
through \ref{Corollary:6} immediately translate over to the
endomorphism ring case.  In particular, we have:

\begin{cor}\label{Corollary:7}
If $M$ has a Dedekind-finite, semi-$\pi$-regular endomorphism
ring, then $M$ has countable exchange.  If, further, the
endomorphism ring is semi-regular, then $M$ has full exchange.
\end{cor}

\section*{\S 7. Final Remarks}\label{Section:FinalRemarks}

In \cite{Ni}, we define what we call {\it finitely complemented}
modules.  These are modules whose direct summands have only
finitely many complement summands.  We showed that a finitely
complemented module with a regular endomorphism ring has full
exchange. We claim that using the methods derived above, one can
remove the condition that $E$ is regular, and replace it with $M$
having finite exchange and $(C_{2})$.

There is another class of modules we can apply these techniques
to; namely, square-free modules.  Suppose that $M$ is a
square-free module with finite exchange.  Mohamed and M\"{u}ller
have shown that $E/J(E)$ is abelian, \cite[Lemmata 11 and
15]{MM1}.  In particular, the element
$\varphi=\sum_{i<\alpha}e_{i,i}$, used in our proof above, is an
idempotent in $E/J(E)$. [Since
$\overline{e}_{i,i}\overline{e}_{j,j}=0$ for $i<j$, and since
idempotents commute in an abelian ring, $\overline{\varphi}$ is a
sum of orthogonal idempotents, and hence is an idempotent.]  Since
idempotents lift modulo $J(E)$ (because $E$ is suitable) we can
lift $\overline{1-\varphi}$ to an idempotent $\widetilde{p}$, as
before. Notice that $\widetilde{\varphi}=\varphi+\widetilde{p}$ is
a unit since $\widetilde{\varphi}$ is congruent to 1 modulo
$J(E)$.  One chooses $p$ as in Theorem~\ref{Theorem:4}. Finally,
if $M$ has $(C_{2})$ we can proceed as in Theorem~\ref{Theorem:6}
to show full exchange for $M$. However, for square-free modules,
the $(C_{2})$ property is equivalent to cohopfianness. So what we
have shown is that a cohopfian, square-free module with finite
exchange has full exchange.

As far as we know, the only classes of modules where it is known
that finite exchange implies countable exchange, but not known if
this further implies full exchange, are square-free modules,
cohopfian modules, and finitely complemented modules.

\end{document}